\documentclass[11pt]{article}
\usepackage{amsfonts,amsmath,latexsym,color,epsfig}
\date{}

    \newcommand{\Dim}{{\mbox{\rm dim }}}
    \newcommand{\Ker}{{\mbox{\rm ker}}}
    \newcommand{\im}{{\mbox{\rm im}}}

\def\qed{\ifhmode\unskip\nobreak\fi\quad\ifmmode\Box\else$\Box$\fi}

\title{On well-connected sets of strings}

\author{
{\sl Peter Frankl}\thanks{R\'enyi Institute, P.O.Box 127 Budapest, 1364 Hungary and MIPT, Moscow, partially supported by the Ministry of Education and Science of the Russian Federation in the framework of MegaGrant 075-15-2019-1926. Email: \texttt{peter.frankl@gmail.com}.}
\and
{\sl J\'{a}nos Pach}\thanks{R\'enyi Institute, P.O.Box 127 Budapest, 1364 Hungary and MIPT, Moscow, partially supported by ERC Advanced Grant GeoScape, NKFIH (Hungarian National Research, Development and Innovation Office) grant K-131529, and by the Ministry of Education and Science of the Russian Federation in the framework of MegaGrant 075-15-2019-1926. Email: \texttt{pach@cims.nyu.edu}}}

\date{}

\begin{document}

\maketitle

\begin{abstract}
Given $n$ sets $X_1,\ldots, X_n$, we call the elements of $S=X_1\times\ldots\times X_n$ {\em strings}. A nonempty set of strings $W\subseteq S$ is said to be {\em well-connected} if for every $v\in W$ and for every $i\, (1\le i\le n)$, there is another element $v'\in W$ which differs from $v$ only in its $i$th coordinate. We prove a conjecture of Yaokun Wu and Yanzhen Xiong by showing that every set of more than $\prod_{i=1}^n|X_i|-\prod_{i=1}^n(|X_i|-1)$ strings has a well-connected subset. This bound is tight.
\end{abstract}

\section{Introduction}\label{intro}

Let $X_1,\ldots, X_n$ be pairwise disjoint sets with $|X_i|=d_i>1$ for $1\le i\le n$. Let $$S=X_1\times\ldots\times X_n=\{(x_1,\ldots,x_n) : x_i\in X_i \mbox{ for every } i\in [n]\}$$
be the set of {\em strings} $x=(x_1,\ldots,x_n)$, where $x_i$ is called the $i$th coordinate of $x$ and $[n]=\{1,\ldots,n\}$.

A subset $W\subseteq S$ is called {\em well-connected} if for every $x\in W$ and for every $i\in [n]$, there is another element $x'\in W$ which differs from $x$ only in its $i$th coordinate. That is, $x'_j\not = x_j$ if and only if $j=i$.

The following statement was conjectured by Yaokun Wu and Yanzhen Xiong~\cite{WX}.

\medskip
\noindent{\bf Theorem 1.} {\em Let $T$ be a subset of $S=X_1\times\ldots\times X_n$ with $|X_i|=d_i>1$ for every $i\in[n]$.
If $$|T|>\prod_{i=1}^nd_i-\prod_{i=1}^n(d_i-1),$$  then $T$ has a nonempty well-connected subset. This bound cannot be improved.}
\medskip

To see the tightness of the theorem, fix an element $y_i$ in each $X_i$ and let $X'_i=X_i\setminus\{y_i\}$. We claim that the set of strings
\begin{equation}\label{eq0}
T_0=(X_1\times\ldots\times X_n)\setminus (X'_1\times\ldots\times X'_n)
\end{equation}
does not have any nonempty well-connected subset. Suppose for contradiction that there is such a subset $W\subseteq T_0$, and let $x=(x_1,\ldots,x_n)$ be an element of $W$ with the minimum number of coordinates $i$ for which $x_i=y_i$ holds. Obviously, this minimum is positive, otherwise $x\not\in T_0$. Pick an integer $k$ with $x_{k}=y_k$. Using the assumption that $W$ is well-connected, we obtain that there exists $x'\in W$ that differs from $x$ only in its $k$th coordinate. However, then $x'$ would have one fewer coordinates with $x_i=y_i$ than $x$ does, contradicting the minimality of $x$.
\smallskip

In the next section, we establish a result somewhat stronger than Theorem 1: we prove that under the conditions of Theorem 1, $T$ also has a subset $W$ such that for every $x\in W$ and $i\in[n],$ the number of elements $x'\in W$ which differ from $x$ only in its $i$th coordinate is {\em odd} (see Theorem 6). In Section~\ref{2ndProof}, we present a self-contained argument which proves this stronger statement.

Shortly after learning about our proof of the conjecture of Wu and Xiong, another proof was found by Chengyang Qian.

\section{Exact sequence of maps}\label{proof}

In this section, we introduce the necessary definitions and terminology, and we apply a basic topological property of simplicial complexes to establish Theorem 1. We will assume throughout, without loss of generality, that the sets $X_i$ are pairwise disjoint.
\smallskip

For every $k\, (0\le k\le n)$, let
     $$S_k=\{ A\subseteq X_1\cup\ldots\cup X_n : |A|=k \mbox{ and } |A\cap X_i|\le 1 \mbox{ for every } i\}.$$
Clearly, we have $|S_n|=|S|=\prod_{i=1}^n|X_i|$. With a slight abuse of notation, we identify $S_n$ with $S$.
The set system $\cup_{k=0}^{n}S_k$ is an \emph{abstract simplicial complex}, that is, for each of its elements $A$, every subset of $A$ also belongs to $\cup_{k=0}^{n}S_k$. This simplicial complex has a geometric realization in $\mathbb{R}^{2n-1}$, where every element $A$ is represented by an $(|A|-1)$-dimensional simplex. (See \cite{B}, part II, Section 9 or \cite{M}, Section 1.5. Note that not all textbooks consider the empty set a $-1$-dimensional simplex, but we do.)
\smallskip 

Assign to each $A\in S_k$ a different symbol $v_A$, and define $V_k$ as the family of all formal sums of these symbols with coefficients $0$ or $1$. Then
     $$V_k=\{ \sum_{A\in S_k}\lambda_A v_A : \lambda_A=0 \mbox{ or } 1 \}$$
can be regarded as a vector space over GF(2) whose dimension satisfies
\begin{equation}\label{eq1}
\Dim V_k=|S_k|=\sum_{1\le j_1< j_2<\ldots< j_k\le n} d_{j_1}d_{j_2}\cdot\ldots\cdot d_{j_k}.
\end{equation}
\smallskip

We use the standard definition of the {\em boundary operations} $\partial_k$. (See~\cite{H}, Section 2.1.) Informally, the boundary of each $(k-1)$-dimensional simplex that corresponds to a member $A\in S_k$ consists of all $(k-2)$-dimensional simplices corresponding to $(k-1)$-element subsets $B\subset A$. This definition naturally extends to any collection (``chain'') of $(k-1)$-dimensional simplices that correspond to members of $S_k$, with multiplicities taken modulo 2.
\medskip

\noindent{\bf Definition 2.} Let $\partial_0:V_0\rightarrow 0$. For every $k\in[n]$ and every $A\in S_k$, let
$$\partial_k(v_A)=\sum_{\substack{B\subset A\\ |B|=k-1}}v_B.$$
Extend this map to a homomorphism $\partial_k:V_k\rightarrow V_{k-1}$ by setting
$$\partial_k(\sum_{A\in S_k}\lambda_A v_A)=\sum_{A\in S_k}\lambda_A\partial_k(v_A),$$
where the sum is taken over GF(2). %If this leads to no confusion, we write $\partial$ for $\partial_k$, for simplicity.
\medskip

Let $\Ker(\partial_k)\subseteq V_k$ and $\im(\partial_k)\subseteq V_{k-1}$ denote the {\em kernel} and the {\em image} of $\partial_k$, respectively.

Our proof is based on the following lemma.

\medskip
\noindent{\bf Lemma 3.} {\em The sequence of homomorphisms
$V_n \xrightarrow{\partial_n} V_{n-1}\xrightarrow{\partial_{n-1}}\ldots\xrightarrow{\partial_1} V_0\xrightarrow{\partial_0} 0$
is an {\em exact sequence}, i.e.,  $\im(\partial_k)=\Ker(\partial_{k-1})$ holds for every $k\in[n]$.}
\medskip

\noindent{\bf Proof.} Before proving the statement, we show that $\im(\partial_k)\subseteq\Ker(\partial_{k-1})$ for every $k\in[n]$. The statement is obviously true for $k=1$.
If $k\ge 2$, then for every $A\in S_k$, we have
$$\partial_{k-1}\partial_{k}v_A=\sum_{\substack{B\subset A\\ |B|=k-1}}\sum_{\substack{C\subset B\\ |C|=k-2}}v_C
=\sum_{\substack{C\subset A\\ |C|=k-2}}2v_C=0.$$
Thus, $\partial_{k-1}\partial_{k}(v)=0$ for every $v\in V_k$, as claimed. In fact, the containment $\im(\partial_k)\subseteq\Ker(\partial_{k-1})$ holds for \emph{every} simplicial complex.
\smallskip

We prove that in our case, all the above containments hold with equality. For every $i\in[n]$, let $K_i$ denote the $0$-dimensional abstract simplicial complex consisting of the $1$-element subsets of $X_i$ and the empty set. Consider now their \emph{join} $K=K_1\ast\ldots\ast K_n$; see~\cite{H}, Chapter 0. By definition, $K$ is the same as the simplicial complex $\cup_{i=0}^nS_i$.
\smallskip

Let $j\ge -1$ be an integer. We need three well known properties of the notion of \emph{$j$-connectedness} of complexes; see Proposition 4.4.3 in \cite{M}.
\begin{enumerate}
\item[(i)]  A complex is $-1$-connected if and only if it contains a nonempty simplex. 
\item[(ii)] If $K_1$ is $a$-connected and $K_2$ is $b$-connected, then their join $K_1\ast K_2$ is $(a+b-2)$-connected.
\item[(iii)] If a complex is $j$-connected, then $\im(\partial_{k})=\Ker(\partial_{k-1})$ holds for every $k, \; 1\le k\le j+2$.
\end{enumerate}

In our case, each $X_i$ is nonempty, hence, by property (i), each $K_i$ is $-1$-connected. By repeated application of (ii), we obtain that  $K=K_1\ast\ldots\ast K_n$ is $(n-2)$-connected. In view of (iii), this implies that $\im(\partial_{k})=\Ker(\partial_{k-1})$ for every $k\in[n]$, as required.  \hfill $\Box$
\medskip

\noindent{\bf Corollary 4.}  {\em For every $k$\; $(0\le k\le n)$, we have\;\;
$\Dim\Ker(\partial_k)=\sum_{i=0}^k(-1)^{k-i}\Dim V_i.$}
\medskip

\noindent{\bf Proof.} By induction on $k$. According to the Rank Nullity Theorem, we have
\begin{equation}\label{eq2}
\Dim V_i=\Dim\Ker(\partial_i)+\Dim\im(\partial_i),
\end{equation}
for every $i\le n$. Since $\Dim V_0=1$ and $\Dim\im(\partial_0)=\Dim 0=0$, the corollary is true for $k=0$.

Assume we have already verified it for some $k<n$. To show that it is also true for $k+1$, we use that
$\Dim\im(\partial_{k+1})=\Dim\Ker(\partial_k)$, by Lemma 3. Plugging this into (\ref{eq2}) with $i=k+1$, we obtain
$$\Dim V_{k+1}=\Dim\Ker(\partial_{k+1})+\Dim\Ker(\partial_k).$$
Hence, using the induction hypothesis, we have
\begin{align*}
\Dim\Ker(\partial_{k+1})&=\Dim V_{k+1}-\Dim\Ker(\partial_k)\\
&=\Dim V_{k+1}-\sum_{i=0}^k(-1)^{k-i}\Dim V_i=\sum_{i=0}^{k+1}(-1)^{k+1-i}\Dim V_i,
\end{align*}
as required. \hfill   $\Box$
\medskip

By (\ref{eq1}), we know the value of $\Dim V_i$ for every $i$. Therefore, Corollary 4 enables us to compute $\Dim\Ker(\partial_n)$ and, hence, $\Dim V_n - \Dim\Ker(\partial_n)$.
\medskip

\noindent{\bf Corollary 5.} {\em We have
$$\Dim V_n - \Dim\Ker(\partial_n)=\prod_{i=1}^nd_i-\prod_{i=1}^n(d_i-1).$$}

\noindent{\bf Proof.} From Corollary 4, we get
$$\Dim V_n-\Dim\Ker(\partial_n)=\sum_{i=0}^{n-1}(-1)^{n-1-i}\Dim V_i.$$
Using (\ref{eq1}) and the fact that $\Dim V_0=1$, this is further equal to
$$\sum_{i=1}^{n-1}(-1)^{n-1-i}\sum_{1\le j_1<j_2<\ldots<j_i\le n}d_{j_1}d_{j_2}\cdots d_{j_i}+(-1)^{n-1}
=\prod_{i=1}^nd_i-\prod_{i=1}^n(d_i-1).\hskip2cm  \hfill \Box$$

\medskip

Now we are in a position to establish the following statement, which is somewhat stronger than Theorem 1.
\medskip

\noindent{\bf Theorem 6.} {\em Let $T$ be a subset of $S=X_1\times\ldots\times X_n$ with $|X_i|=d_i>1$ for every $i\in[n]$.
If $$|T|>\prod_{i=1}^nd_i-\prod_{i=1}^n(d_i-1),$$
then there is a nonempty subset $W\subseteq T$ with the property that for every $x\in W$ and $i\in[n],$ the number of elements $x'\in W$ which differ from $x$ only in their $i$th coordinate is {\em odd}. This bound cannot be improved.}
\medskip

\noindent{\bf Proof.} The tightness of the bound follows from the tightness of Theorem 1 shown at the end of the Introduction.

Let $T$ be a system of strings of length $n$ satisfying the conditions of the theorem. Using the notation introduced at the beginning of this section, let
$$V(T)=\{ \sum_{A\in T}\lambda_A v_A : \lambda_A=0 \mbox{ or } 1 \}.$$
Then $V(T)$ can be regarded as a linear subspace of $V_n$ with $\Dim V(T)=|T|$. Comparing the size of $T$ with the value of $\Dim V_n - \Dim\Ker(\partial_n)$ given by Corollary 5, we obtain that there is a nonzero vector
$v=\sum_{A\in T}\lambda_A v_A$ that belongs to $V(T)\cap\Ker(\partial_n)$. Let $W=\{ A\in T : \lambda_A=1 \}$. Then we have
$$0=\partial_n(v)=\sum_{A\in W}\partial_n(v_A)=\sum_{A\in W}\sum_{\substack{B\subset A\\ |B|=n-1}}v_B
=\sum_{\substack{B\subset [n]\\ |B|=n-1}}|\{A\in W : A\supseteq B\}|v_B.$$
Thus, for each $B$, the coefficient of $v_B$ is {\em even}. This means that the set of strings $W\subset T$ meets the requirements of the theorem. \hfill  $\Box$

\section{Direct proof of Theorem 6}\label{2ndProof}

In this section, we prove Corollary 5 and, hence, Theorem 6 directly, without using Lemma 3.
\smallskip

As in the Introduction, fix an element $y_i\in X_i$ and let $X'_i=X_i\setminus\{y_i\}$, for every $i\in[n]$. Defining $T_0$ as in (\ref{eq0}), we have that $|T_0|=\prod_{i=1}^nd_i-\prod_{i=1}^n(d_i-1)$.

Suppose that $|T|>|T_0|$. To prove Corollary 5, it is sufficient to show that there exists a nonzero vector $v=\sum_{A\in T}\lambda_A v_A$ with suitable coefficients $\lambda_A\in\{0,1\}$ such that $v\in\Ker(\partial_n)$, {\em i.e.}, we have $\partial_nv=\sum_{A\in T}\lambda_A(\partial_n v_A)=0$. Thus, it is enough to establish the following statement.
\medskip

\noindent{\bf Lemma 7.} {\em Let $T$ be a subset of $S=X_1\times\ldots\times X_n$ with $|X_i|>1$ for every $i\in[n]$.

If $|T|>|T_0|$, then the set of vectors $\{\partial_{n}v_A : A\in T\}$ is linearly dependent over {\rm GF(2)}.}
\medskip

\noindent{\bf Proof.}  First, we show that the set of vectors $\{\partial_{n}v_A : A\in T_0\}$ is linearly independent. Suppose, for a contradiction, that there is a nonempty subset $W\subset T_0$ such that
$\sum_{A\in W}\partial_{n}v_A=0$. Pick an element $A=\{x_1,\ldots,x_n\}$ of $W$ for which the number of coordinates $i$ with $x_i=y_i$ is as small as possible. By the definition of $T_0$, there is at least one such coordinate $x_k=y_k$. In view of Definition 2, one of the terms of the formal sum $\partial_n v_A$ is
$v_{B}$ with $B=A\setminus\{y_k\}$, and this term cannot be canceled out by a term of $\partial_n v_{A'}$ for any other $A'\in W$, because in this case $A'$ would have fewer coordinates that are equal to some $y_i$ than $A$ does. Hence,
$\sum_{A\in W}\partial_{n}v_A\not=0$, contradicting our assumption.
\smallskip

It remains to prove that $\{\partial_{n}v_A : A\in T_0\}$ is a \emph{base} of $\im(\partial_n)$, that is,  there exists no set of strings $T\supset T_0$ with $|T|>|T_0|$ such that the set of vectors $\{\partial_{n}v_A : A\in T\}$ is linearly independent over GF(2).

To see this, consider any string $C=\{z_1,\ldots,z_n\}\in S\setminus T_0$. Since $C\not\in T_0$, we have $z_i\not=y_i$ for every $i$. Define $T(C)$ as the set of all strings $A=\{x_1,\ldots,x_n\}\in S$ whose every coordinate $x_i$ is either $y_i$ or $z_i$. Then we have $\sum_{A\in T(C)}\partial_{n}v_A=0$. As we have $T(C)\subseteq T_0 \cup\{C\}$, this means that the set of vectors $\{\partial_{n}v_A : A\in T_0\cup\{C\}\}$ is linearly dependent over GF(2). This completes the proof of the lemma and, hence, of Theorem 6. \hfill $\Box$

\medskip

\noindent\textbf{Acknowledgement.} We thank G\'abor Tardos and an anonymous referee for several helpful suggestions.

\end{document}